\documentclass{ifacconf}
\usepackage{amssymb}
\usepackage{amsmath}
\usepackage{color}
\usepackage{graphicx}
\usepackage{caption}
\usepackage{subcaption}
\usepackage{tikz}
\usetikzlibrary{positioning, arrows.meta}
\usetikzlibrary{arrows.meta,positioning,calc,shapes.geometric,fit,backgrounds,matrix}
\usetikzlibrary{decorations.pathreplacing}
\newtheorem{example}{Example}

\usepackage{graphicx}      
\usepackage{natbib}        
\usepackage{float}
\usepackage{caption} 
\usepackage[normalem]{ulem}
\usepackage{xcolor}

\newcommand{\SY}{\color{red}}

\begin{document}
\raggedbottom
\begin{frontmatter}

\title{Physics Informed Neural Networks for Nonlinear Delay Differential Equations} 


\author[First]{Stone Yao} 
\author[Second]{Vipin Kumar} 
\author[Third]{Roberto Guglielmi}

\address[First]{Department of Computer Science, University of Waterloo,\\ Waterloo, ON, 
Canada (e-mail: stone.yao@uwaterloo.ca).}
\address[Second]{Department of Mathematics \& Computing, Dr B. R. Ambedkar National Institute of Technology Jalandhar, Punjab, India\\ (e-mail: vipinkumar@nitj.ac.in)}
\address[Third]{Department of Applied Mathematics, University of Waterloo, Waterloo, ON, Canada (e-mail: roberto.guglielmi@uwaterloo.ca)}

\begin{abstract}                
In this paper we propose a novel physics-informed neural network framework for solving general first-order delay differential equations. Our approach combines a differentiable history switch, a trial-solution formulation that explicitly enforces history constraints, and a segmented collocation strategy to stabilize gradient propagation across large temporal domains. The method enables a 
scalable and physics-consistent approximation of delay differential equation solutions while maintaining continuity across subintervals. Numerical experiments demonstrate the effectiveness of the proposed method.
\end{abstract}


\end{frontmatter}

\section{Introduction}
Differential equations arise in many physical applications of science and engineering, such as population dynamics, mechanical systems, electrical circuits, and biological processes. In many real-world phenomena, the dynamics of a system at a given time depend not only on its current state but also on its past history—introducing delay into the dynamics. Such systems are modelled by delay differential equations (DDEs), which naturally occur in diverse areas including population dynamics, epidemiology, neural networks, control theory, and economics, see \cite{book-delay-1}. In contrast to ordinary differential equations (ODEs), the presence of delay introduces infinite-dimensional dynamics, leading to richer but more complex behaviours such as oscillations, bifurcations, and stability switches, see \cite{effect-delay-1,effect-delay-2,effect-delay-3}. 

An accurate method for approximating the solutions of DDEs is therefore essential for understanding and predicting the evolution of time-delay systems, and shall take these challenges into account. Indeed, analytical methods for DDEs are typically limited to specific forms or linearized systems, often relying on characteristic equations or Laplace transforms, see \cite{book-ODEs}. For most nonlinear and non-autonomous problems, closed-form solutions are unavailable, motivating the development of numerical methods. Classical numerical solvers, such as the method of steps, Runge–Kutta-based schemes, or collocation techniques (e.g., MATLAB’s \texttt{dde23}), approximate the solution by discretizing time and propagating it iteratively using previously computed delay values, see \cite{book-numerical-odes-1,book-numerical-odes-2}. While effective, these methods can become computationally expensive or unstable for long-time horizons, high-dimensional systems, or when the delay term is state-dependent or distributed.

Early neural approaches for differential equations date back to the 1990s, where feed-forward networks were trained to minimize residual errors of ODEs and PDEs \cite{Lagaris, Meade, Dissanayake}. Subsequent works have extended these ideas to various equation types, including stiff ODEs \cite{stiff-ode}, fractional systems \cite{frac}, and delay equations \cite{delay-ode-1,delay-ode-2,delay-Noor}. For DDEs, some recent studies have used recurrent neural networks or hybrid architectures combining history functions with data-driven delay modelling \cite{delay-ode-3,delay-ode-4}. However, most of these methods are problem-specific, require large training datasets, or lack rigorous enforcement of physical laws.

The Physics-Informed Neural Network (PINN) framework, introduced by~\cite{Raissi}, provides a general methodology for solving differential equations by incorporating the governing physical laws into the neural network’s loss function {while retaining the benefit of automatic differentiation in the training phase}. PINNs have been successfully applied to a wide range of problems, including forward and inverse modelling of ODEs and PDEs \cite{Pinns-app-1,Pinns-app-2}, fractional systems \cite{Pinns-app-frac-1,Pinns-app-frac-2}, stochastic differential equations \cite{Pinns-app-stochastic-1}, control systems \cite{Pinns-app-control} and many more. 
While PINNs have been extensively studied for ODEs and PDEs, their application to DDEs remains relatively limited. A few recent efforts have incorporated delay terms by augmenting the network inputs with historical states or by using recurrent formulations \cite{Pinns-delay-1, Pinns-delay-2}. 

To the best of our knowledge, this is the first contribution to develop a comprehensive PINN-based framework designed to solve general first-order DDEs with arbitrary history functions and delay terms in a unified, stable, and scalable manner. Most existing implementations are case-dependent and fail to generalize across different delay magnitudes or dynamic regimes.
In this work, we propose a novel PINN framework for solving general DDEs by integrating a differentiable history switch, a trial-solution formulation that enforces history constraints explicitly, and a segmented collocation strategy to stabilize gradients across large temporal domains. Our approach ensures continuity between sub-intervals without increasing model complexity and avoids unstable implicit gradients through delayed dependencies by stopping gradient flow in delayed branches. The resulting method provides a mesh-free, scalable, and physics-consistent neural approximation that achieves high accuracy.


\section{Problem Definition and Methodology}\label{sec:ProbForm}
In this section, we introduce the key components of the proposed segmented Trial-PINN framework. We begin with the problem setup and the trial-solution construction, then describe domain segmentation and collocation strategies, lastly, the neural approximator, the differentiable history switch used for handling delays, and the terms in the physics-based loss, as well as the optimization procedure.
\begin{figure}[b]
\centering
\includegraphics[width=\linewidth]{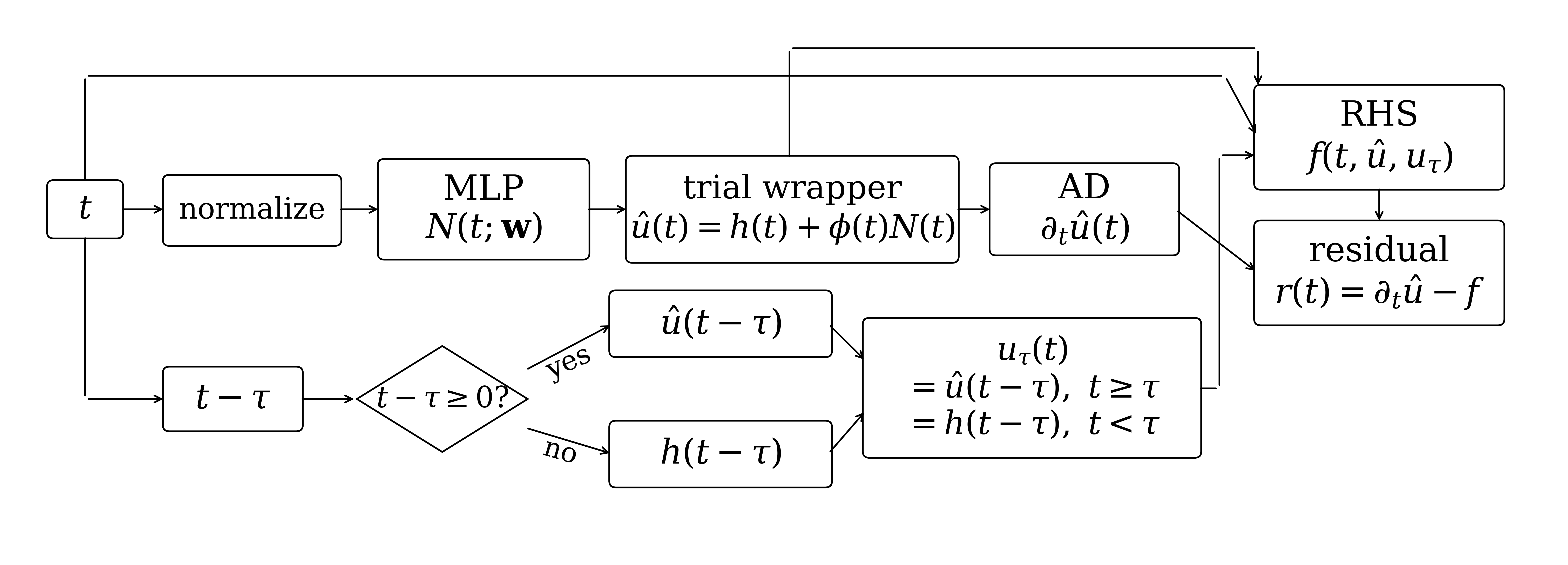}
\caption{Trial-PINN dataflow and residual construction for general DDEs as in~\eqref{eq:generalDDE}.}
\label{fig:figure1}
\end{figure}

$\textbf{Problem Setup}$ \newline
We consider first-order nonlinear DDEs of the form:
\begin{equation}\label{eq:generalDDE}
    \begin{array}{ll}
        u'(t) = f(t,u(t),u(t-\tau)), \quad &t\in(0,T], \\[1ex]
        u(t) = h(t), \quad &t\in[-\tau,0],
    \end{array}
\end{equation}
where $u:[-\tau,T] \to \mathbb{R}^n$ is the state, {$\tau>0$ is a constant delay, $f$ is the nonlinear dynamics depending on the delayed state $u(\cdot -\tau)$ besides possibly the independent variable $t$ and the current state $u(t)$,} $h(t)$ is the prescribed history function, and $T>0$ is the final time.
To approximate the solution, we construct a trial solution that enforces the history exactly:
\begin{equation}\label{eq:trial-sol}
\begin{array}{ll}
\!\!\!   \hat{u}(t,w)=h(t)+\phi(t)N(t,w), \quad & t\in [-\tau,T], w\in \mathbb{R}^m,\\[1ex]
    \phi(t):=\frac{\operatorname{max}(0,t)}{T}, \quad & t\in [-\tau,T],
\end{array}
\end{equation}
where $N(t, w)$ denotes a fully connected neural network with weights $w$. The function $\phi(t)$ satisfies $\phi(t) = 0$ on the history interval, ensuring that $\hat{u}(t)=h(t)$ holds for $t\leq0$ and for $t>0$, $\phi(t)$ gradually activates the neural network component.
The overall data-flow of our architecture consists of: normalization, neural prediction, trial wrapping, delay evaluation, automatic differentiation, and residual construction. These components are illustrated in Figure $\ref{fig:figure1}$, and we will describe some in more detail in the following paragraphs.

\begin{figure}[H]
\centering
\resizebox{\linewidth}{!}{%
\begin{tikzpicture}[
  node distance = 8mm and 10mm,
  >={Latex},
  block/.style      = {rectangle, draw, rounded corners, align=center,
                       inner sep=3pt, minimum height=6mm, font=\footnotesize,
                       text width=32mm, fill=white},
  smallblock/.style = {rectangle, draw, rounded corners, align=center,
                       inner sep=2pt, minimum height=5mm, font=\scriptsize,
                       text width=40mm, fill=white},
  tinyblock/.style  = {rectangle, draw, rounded corners, align=center,
                       inner sep=2pt, minimum height=5mm, font=\scriptsize,
                       text width=12mm, fill=white},
  decision/.style   = {diamond, draw, aspect=2, align=center, inner sep=1pt,
                       font=\scriptsize, fill=white},
  arrow/.style      = {->, line width=0.5pt}
]
  \draw (0,0) -- (10,0);
  \foreach \x/\lab in {0/$0$, 1.5/$t_1$, 3/$t_2$, 4.5/$t_3$, 8.5/$t_{K-1}$, 10/$T$}
    {\draw (\x,0.08) -- (\x,-0.08) node[below] {\lab};}
  \node[above] at (0.75,0.2) {$I_1$};
  \node[below] at (0.75,-0.3) {$C_1$};
  \draw[decorate,decoration={brace,amplitude=5pt,mirror}] (0,-0.15) -- (1.5,-0.15) node[midway,yshift=-0.6cm] {};
  \node[above] at (2.25,0.2) {$I_2$};
  \node[below] at (2.25,-0.3) {$C_2$};
  \draw[decorate,decoration={brace,amplitude=5pt,mirror}] (1.5,-0.15) -- (3,-0.15) node[midway,yshift=-0.6cm] {};
  \node[above] at (3.75,0.2) {$I_3$};
  \node[below] at (3.7,-0.3) {$C_3$};
  \draw[decorate,decoration={brace,amplitude=5pt,mirror}] (3,-0.15) -- (4.5,-0.15) node[midway,yshift=-0.6cm] {};
  \node[above] at (9.25,0.2) {$I_K$};
  \node[below] at (9.25,-0.3) {$C_K$};
  \draw[decorate,decoration={brace,amplitude=5pt,mirror}]
  (8.5,-0.15) -- (10,-0.14) node[midway,yshift=-0.6cm] {}; \node[above] at (6.5,0.2) {$\dots$};
  \foreach \x in {0.1,0.2,0.3,0.4,0.5,0.6,0.7,0.8,0.9,1,1.1,1.2,1.3,1.4,1.5} \fill (\x,0) circle (1pt);
  \foreach \x in {1.6,1.7,1.8,1.9,2.0,2.1,2.2,2.3,2.4,2.5,2.6,2.7,2.8,2.9,3.0} \fill (\x,0) circle (1pt);
  \foreach \x in {3.1,3.2,3.3,3.4,3.5,3.6,3.7,3.8,3.9,4.0,4.1,4.2,4.3,4.4,4.5} \fill (\x,0) circle (1pt);
  \foreach \x in {8.6,8.7,8.8,8.9,9.0,9.1,9.2,9.3,9.4,9.5,9.6,9.7,9.8,9.9,10} \fill (\x,0) circle (1pt);
\end{tikzpicture}
}
\caption{Domain decomposition of $[0,T]$ into $K$ subintervals $I_k$ with independent collocation sets ${C}_k := \{c_{j}\}_{j=1}^{n_k}\subset I_k$ for each $k=1,\ldots,K$, with $n_k$ defined according to the paragraph {\bf Collocation Strategies}.}\label{fig:Fig2}
\end{figure}
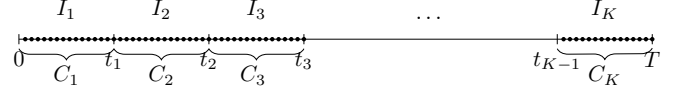

$\textbf{Domain Segmentation}$ \newline
To improve stability and scalability on domains with a large final time $T$, the time horizon interval $[0,T]$ is partitioned into $K$ sub-intervals $I_k$, each with its own collocation set $C_k \subset I_k$, as shown in Figure~\ref{fig:Fig2}.
Segmentation reduces error accumulation and stiffness over long horizons by training one small neural network per subinterval, allowing each region to learn its dynamics independently. After training, these subnetworks are assembled into a single piecewise global model, which improves stability and accuracy over long domains. 

$\textbf{Collocation Strategies}$ \newline
In order to construct the physics-informed loss, we must choose $n_k$ collocation time instances inside each sub-interval $I_k$, $k=1,\ldots,K$. Our framework supports three sampling strategies, each serving a different purpose:
\begin{itemize}
    \item  {Uniform:} Sampled points are evenly distributed across the interval; suitable when the residual behaves smoothly.
    \item  {Non-uniform:} The sampling applies
    \begin{equation}\label{eq:non-unifSam}
    u \mapsto 1 -(1-u)^p, \quad p>1,
    \end{equation}
    to bias samples toward later times, improving resolution in stiff or transition regions. Here, $u \in [0,1]$ is a uniform random variable that is mapped into a time point t.
    \item  {Segmented builder:} A builder returns $K$ independent sets ${\{C_k}\}$, each sampled uniformly or via the above non-uniform rule. This offers flexible per-segment densities and is essential for large-domain problems.
\end{itemize}

\begin{figure}[H]
\centering
\includegraphics[width=\linewidth]{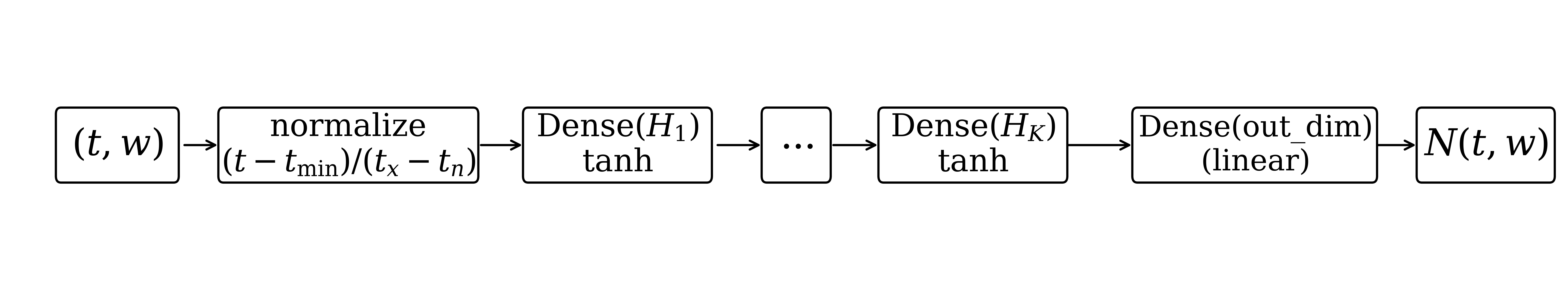}
\caption{Neural approximator $N(t,w)$ used in~\eqref{eq:trial-sol}: input $t$ is normalized, processed by tanh-activated dense layers, and mapped to a linear output head.}
\label{fig:figure3}
\end{figure}

$\textbf{Neural Approximator}$\\ 
The network $N(t;w)$ is a fully connected MLP with tanh activations. The input time $t$ is first normalized to $[0,1]$, passed through several dense layers with fixed width, and mapped to an output through a linear layer (Figure $\ref{fig:figure3}$). The network output $N(t;w)$ is then combined with the scaling function $\phi(t)$ and the history function through the trial wrapper in equation~$\eqref{eq:trial-sol}$ to produce the trial solution~$\hat{u}(t)$, as illustrated in the figure below.

$\textbf{Delay Handling via a Differentiable History Switch}$ \newline
At each collocation point $t_i \in C_k$, evaluation of $u(t_i-\tau)$ is required. We compute this using a differentiable history switch implemented with a TensorFlow library (denoted by \texttt{tf} in the rest of the paper).

\begin{figure}[H]
\centering
\includegraphics[width=\linewidth]{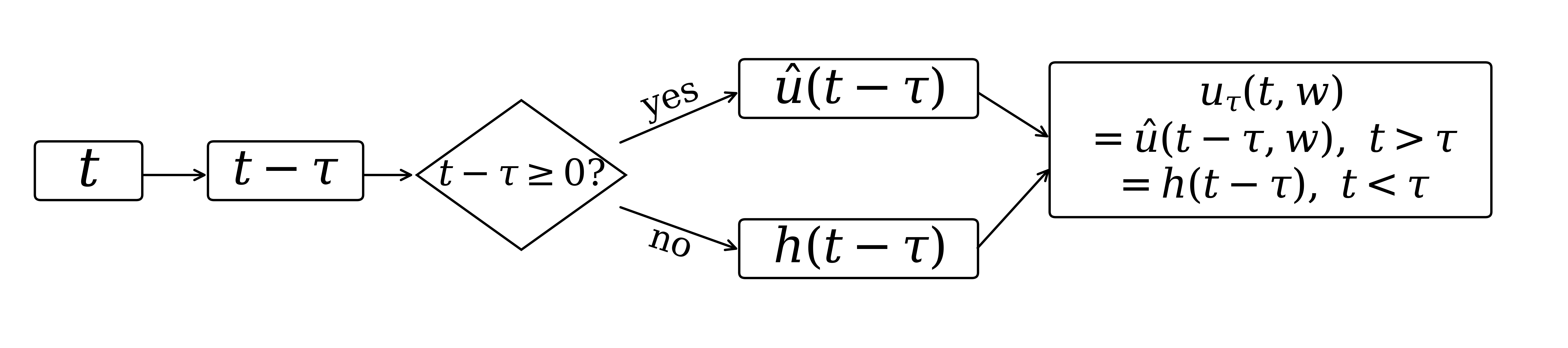}
\caption{Differentiable history switch used for delayed evaluations.}
\label{fig:figure4}
\end{figure}

To prevent long gradient-dependency chains across segments, delayed values from earlier segments are passed through { \texttt{tf.stop\_gradient}}. This mechanism corresponds to the delay branch shown in Figure $\ref{fig:figure1}$ (bottom branch).

$\textbf{Physics Loss Functions}$ \newline
Let $C= \cup_{k=1}^K C_k$ be the full collocation set on the domain $[0,T]$. For any $t\in C$, the physics residual is defined as
\begin{equation*}
    r_k(t, w)=\partial_t \hat{u}(t, w) - f(t,\hat{u}(t, w), u_{\tau}(t,w))
\end{equation*}
where only the temporal derivative of the trial solution $\partial_t \hat{u}(t, w)$ is computed via automatic differentiation. $f$ itself is evaluated but not differentiated with respect to time. Here, $\tau$ denotes the fixed delay of the DDE. The history loss operates on our collocation set of the history domain $\mathcal{H}$ and the true history value at $t_i$, denoted by $y_i$. 
\begin{equation*}
    \mathcal{L}_{hist}(w) = \frac{1}{|\mathcal{H}|} \sum_{t_i \in \mathcal{H}}{(\hat{u}(t_i, w)-y_i)^2}.
\end{equation*}
Although the trial solution already enforces $\hat{u}(t) = h(t)$ for $t \leq 0$, this term provides additional numerical stability near the boundary points and suppresses optimizer drift that may cause deviation from the prescribed history. We also include an optional $\mathcal{L}_{bc}$ loss at boundary points $t_k^*$ between segments: $\mathcal{L}_{bc}(w) = \sum_{k=2}^{K}|\hat{u}_k(t_k^*,w)-\hat{u}_{k-1}(t_k^*,w)|^2$. This soft constraint mitigates interface discontinuities across subdomains. \newline
The total training loss is then:
\begin{multline}\label{eq:lossfunc}
 \mathcal{L}(w) =\lambda_{dde}\mathcal{L}_{dde}(w) + \lambda_{hist}\mathcal{L}_{hist}(w) + \lambda_{bc}\mathcal{L}_{bc}(w)
\end{multline}
We emphasize that setting $K=1$ recovers the standard non-segmented PINN loss, as noted earlier in the domain-segmentation discussion.

$\textbf{Optimization and Implementation Details}$ \newline
Training uses the Adam optimizer with exponential decay. Gradient clipping is used to improve stability, {and the entire loop is compiled using {\texttt{tf.function}}}. Adjustable hyperparameters consist of the number of segments $K$ (for segmented training), the number of epochs, the learning rate, the collocation density and the bias exponent $p$ in~\eqref{eq:non-unifSam} when using non-uniform sampling and weights, $\lambda_{dde}, \lambda_{hist}, \lambda_{bc}$ from~\eqref{eq:lossfunc}.

All first-order and delayed derivatives are computed using automatic differentiation, with delay queries handled through the differentiable history switch introduced earlier (Fig. \ref{fig:figure4}). Under domain segmentation, each subdomain $C_k$ is associated with its own neural submodel, which is trained independently using only the collocation points inside that segment. After training, the final predictor is assembled by concatenating the submodels into a single piecewise predictor, ensuring a continuous piecewise trial solution over $[0,T]$. During training, we record the physics residuals and the history loss to monitor convergence.

{\bf Evaluation}\\
To assess model accuracy, we compare the predicted values $\hat{y}_i$ against reference values $y_i$ collected at $N$ uniformly sampled evaluation points on the domain of interest. These reference values come either from an analytic solution (when available) or from a high-accuracy numerical solution obtained using MATLAB’s \texttt{dde23}. The two error metrics used are the Mean Absolute Error (MAE) and Root Mean Squared Error (RMSE):
\begin{equation*}
    \text{MAE} = \frac{1}{N} \sum^{N}_{i=1} |y_i- \hat y_i |,\ \ 
    \text{RMSE} = \sqrt{\frac{1}{N} \sum^{N}_{i=1}(y_i- \hat y_i)^2}.
\end{equation*}
The MAE quantifies the average point-wise deviation between the predicted and reference solutions, while the RMSE penalizes larger errors more strongly and therefore reflects overall stability of the approximation. To isolate the effect of segmentation, all experiments fix the neural-network architecture and optimizer settings, varying only the number of segments $K$ and the corresponding collocation-point allocation.
\section{Numerical Results}
To evaluate the Segmented Trial-PINN illustrated in Section~\ref{sec:ProbForm}, we test the method by comparing it against the solutions for several representative DDEs, where the solutions are either analytic or high-accuracy numerical references from \texttt{dde23} with tolerances $10^{-8}$ (relative) and $10^{-10}$ (absolute). In all the examples below, we also compare our method to a non-segmented Trial-PINN baseline. The experiments are performed using Python 3 and TensorFlow (CPU) on an Intel Core i7-12700H. {Unless specified otherwise,} the neural network is a fully connected MLP with four 10-neuron hidden layers, tanh activations, Xavier initialization, and Adam (learning rate 0.02 with 30\% decay every 500 epochs over 2500 epochs) with gradient clipping at 100. Any deviations from this baseline, such as larger networks or extended training, are reported in the corresponding examples. Collocation points follow the experiment’s sampling choice (uniform, non-uniform, or segmented), which independently determines interior point selection and domain segmentation.


\begin{example}\label{ex:harmonic} 
We consider the harmonic DDE considered in~\cite{Pinns-delay-2}:
\begin{equation} \label{DDE-1}
\begin{aligned}
    y'(t)&=y(t)+y(t-{\SY {\tau}})-3\cos(t)-5\sin(t), \\[1ex]
y(t)&=3\cos(t)-5\sin(t) \quad \text{for $t \in$ } [-\pi,0],
\end{aligned}
\end{equation}
with delay $\tau = \pi$ and time horizon $T = 20$. This DDE admits the exact solution $y(t)=3 \cos t - 5\sin t$, { $t\in [0,20]$}. The domain $[0, 20]$ is divided into four sub-intervals, each trained independently, using 1000 interior points and 600 history points. The loss is
 $   \mathcal{L}_{total} = \mathcal{L}_{hist} + 2.5\mathcal{L}_{DDE} + 2\mathcal{L}_{bc}.$

Figure~$\ref{fig:example1}$ compares the segmented Trial-PINN, the non-segmented Trial-PINN, and the analytic solution. Over the full domain $[0,20]$, both models accurately reproduce the periodic oscillatory behavior. The non-segmented model maintains amplitude and phase most uniformly, while the segmented model exhibits a small localized deviation near $t \approx 5$ before quickly realigning. Although this slight mismatch arises from subdomain transitions interacting with the delayed term, it remains transient and does not affect global accuracy. Importantly, this example is not one where segmentation is expected to outperform the standard PINN, since the harmonic DDE is smooth, periodic, and relatively well-conditioned. However, the benefits of segmentation become clearer in more complex nonlinear settings, such as the two examples that follow, where sharper transitions, stronger delay feedback, or stiffness create local training difficulties that segmentation is specifically designed to mitigate.
\begin{figure}[h]
  \centering
  \includegraphics[width=\columnwidth]{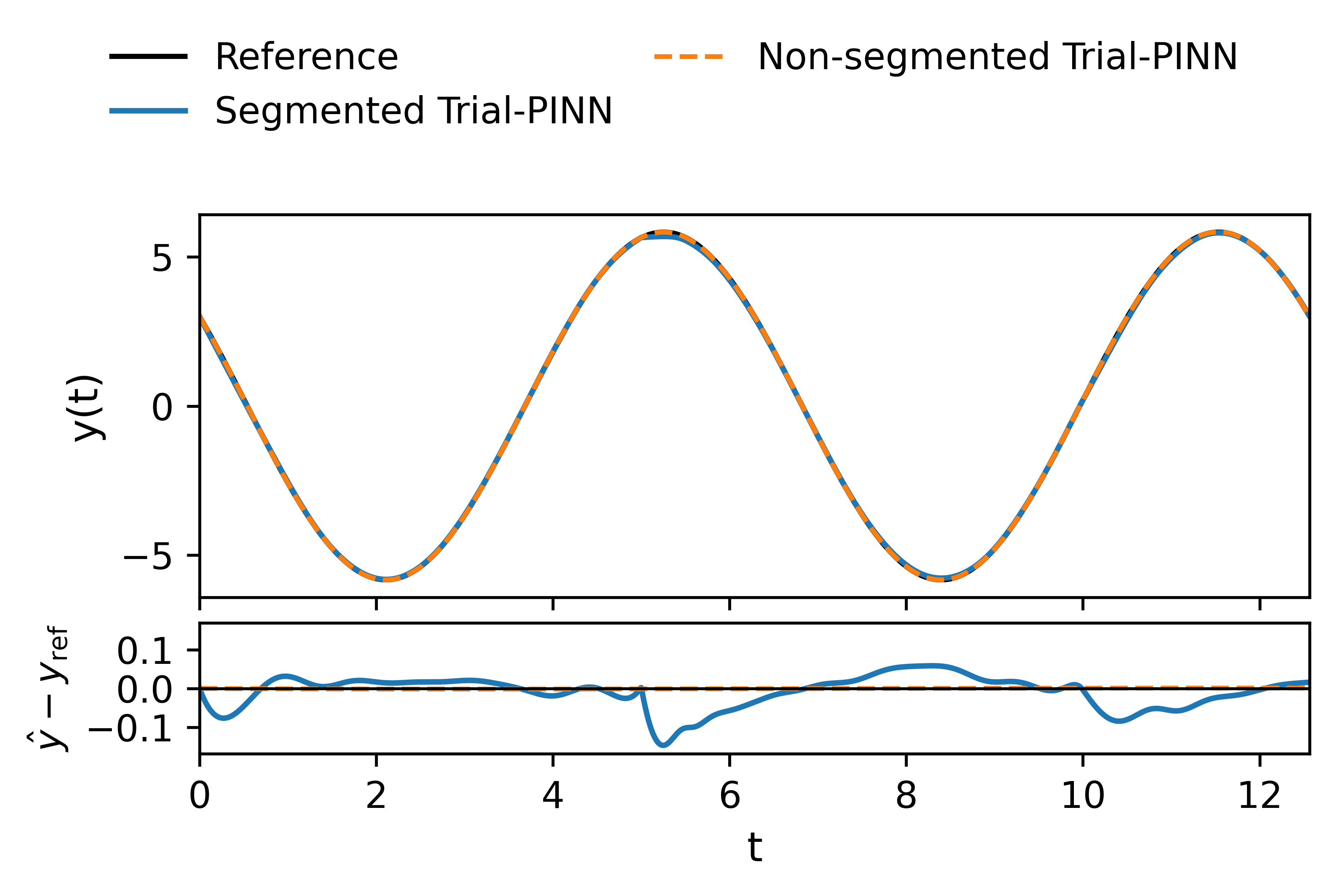}
  \caption{Comparison of the analytic and predicted solutions for the harmonic DDE \eqref{DDE-1}.}
  \label{fig:example1}
\end{figure}
While segmentation offers no advantage here, it provides a scalable and robust framework for harder DDEs. \newline
Table~\ref{tab:training_results_1} compares RMSE, MAE, and training time for segmented and non-segmented Trial-PINN methods. Both models achieve low errors; the non-segmented version attains the smallest RMSE/MAE due to its continuous global representation. The segmented model incurs roughly double the training time because each sub-network is optimized independently, but remains practical for parallelization or adaptive refinement workflows.
\begin{table}[h]
\centering
  \caption{Quantitative comparison on the harmonic DDE \eqref{DDE-1} benchmark.
  }
  \label{tab:training_results_1}
  \begin{tabular}{p{2.85cm} p{1.5cm} p{1.5cm} p{1.2cm}}
    \hline
    Model & RMSE vs Analytic & MAE vs Analytic & Training Time (s) \\
    \hline
    Seg. Trial-PINN & 0.04213 & 0.00177 & 59.10 \\
    Non-seg. Trial-PINN & $1.08 \times 10^{-3}$ & $1.17 \times 10^{-6}$ & 22.23 \\
    \hline
  \end{tabular}
\end{table}

  
\end{example}

\begin{example} \label{ex:logistic}
We examine the classical logistic DDE considered in~\cite{baker-rost-delayed-logistic},
\begin{equation}\label{DDE-2}
\begin{aligned}
    y{'}(t)&=y(t)(1-y(t-\tau)), \quad t\in (0,10),\\
    y(t)&=5, \quad t \in [-2,0],
\end{aligned}
\end{equation}
with delay $\tau = 2$. This nonlinear delayed-feedback system is a standard benchmark due to its rich transient behavior. The domain $[0,10]$ is divided into $K=10$ segments, each trained with $N_f=1000$ interior collocation points per segment and $N_h=2000$ history samples. The total loss takes the weighted form
 $   \mathcal{L}_{total} = \mathcal{L}_{hist} + 2.5\mathcal{L}_{DDE} + 2\mathcal{L}_{bc}.$
To illustrate the benefit of segmentation, we compare our model with three non-segmented Trial-PINN baselines of increasing capacity: PINN-1 has the same architecture as the segmented model, PINN-2 uses a $5\times15$ network with 5000 epochs and PINN-3 uses a $5\times20$ network with 10000 epochs.  \newline
Figure $\ref{fig:example2}$ shows that all models capture the characteristic decaying–oscillatory dynamics. The segmented Trial-PINN {produces a smoother early transient for $t<2$} and remains closely aligned with the reference throughout the horizon $[0,10]$. The smallest non-segmented model struggles with steep gradients, while deeper variants improve accuracy but require substantially more training.
\begin{figure}[H]
  \centering
  \includegraphics[width=\columnwidth]{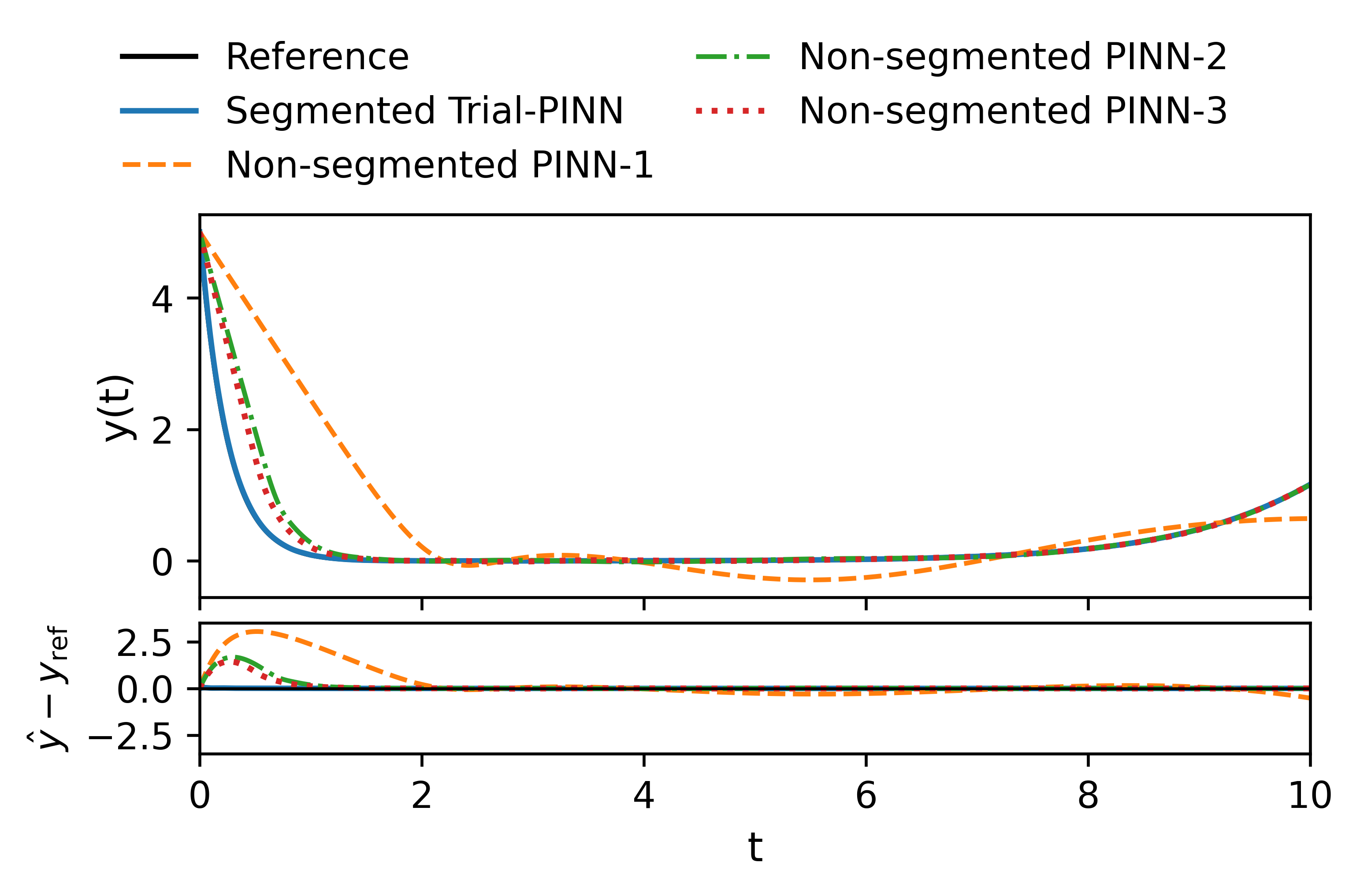}
  \caption{Comparison of segmented and non-segmented models for DDE \eqref{DDE-2}. 
 }
  \label{fig:example2}
\end{figure}
Table~\ref{tab:training_results_3} shows that the segmented model achieves the lowest overall error with moderate computational overhead. Larger non-segmented networks improve the accuracy over smaller non-segmented models, although they remain far less accurate than the segmented Trial-PINN and become significantly more expensive in the training phase. Despite increased training cost, the segmented approach yields over 50× lower RMSE than any non-segmented baseline.
\begin{table}[h]
  \centering 
  \caption{Quantitative comparison of different models on the logistic DDE \eqref{DDE-2} 
  }
\label{tab:training_results_3}
  \begin{tabular}{p{3.1cm} p{1.5cm} p{1.3cm}  p{1.2cm}}
    \hline
    Model & {\scriptsize RMSE ($×10^{-2}$)} & {\scriptsize MAE ($×10^{-5}$)} & {\scriptsize Training Time (s)}\\
    \hline
    $\text{\scriptsize Seg. Trial-PINN}$ & 0.485 & 2.36 & 165.57 \\
    $\text{\scriptsize Non-seg. Trial-PINN-1}$ & 73.78 & 54437.13 & 40.19 \\
    $\text{\scriptsize Non-seg. Trial-PINN-2}$ & 34.64 & 11999.47 & 105.99 \\
    $\text{\scriptsize Non-seg. Trial-PINN-3}$ & 27.94 & 7806.80 & 329.01 \\
    \hline
  \end{tabular}
\end{table}

To further investigate the effectiveness of the segmented method on different delay values, we evaluated the logistic DDE example with two additional delay values, $\tau = 1$ and $\tau = 4$. Figure $\ref{fig:example4}$ shows the comparison between the segmented and non-segmented Trial-PINN predictions against the MATLAB reference solutions. This additional experiment illustrates how the method behaves with shorter or longer delay values, and provides additional evidence of the robustness of the segmented framework.
\begin{figure}[tb]
  \centering
  \includegraphics[width=\columnwidth]{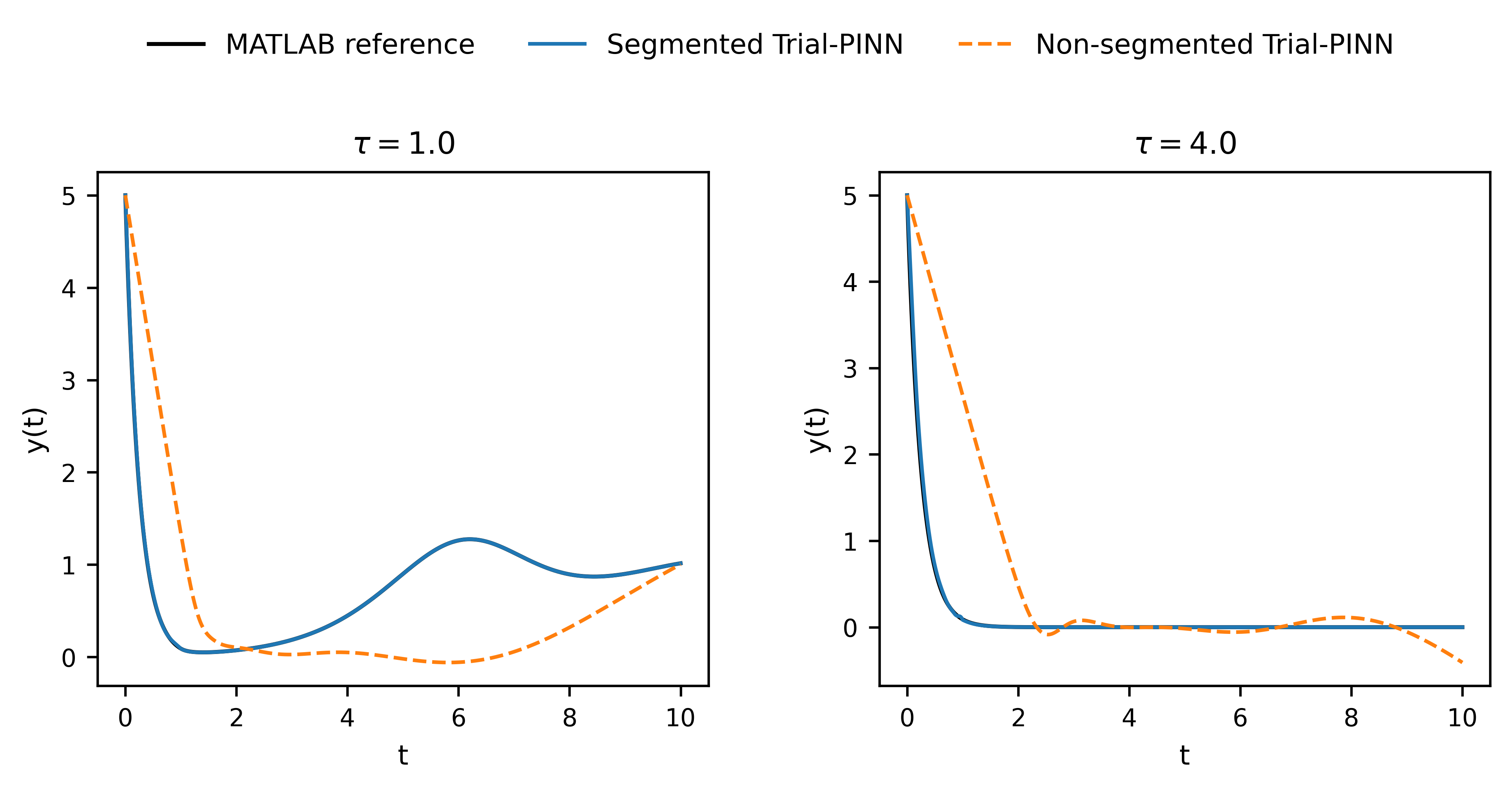}
  \caption{Logistic DDE results for $\tau$ = 1 and $\tau$ = 4, comparing segmented and non-segmented Trial-PINNs with the same hyperparameter settings.
  }
  \label{fig:example4}
\end{figure}
\end{example}

\begin{example} \label{ex:coupled}
To demonstrate scalability, we consider the coupled nonlinear DDE system on the horizon $t\in[0,20]$
\begin{equation} \label{DDE-3}
    \begin{split}
        u'(t) &= -\alpha u(t) + \beta v(t-\tau) + \gamma \sin(t), \\
        v'(t) &= \delta u(t-\tau)-\eta v(t) + \kappa \cos(t),
    \end{split}
\end{equation}
with history functions
 $u(t) = \sin(t), \quad v(t) = \cos(t), \quad t \in [-\tau,0], \quad \tau=1.5,$
and parameters $\alpha = 0.8, \beta = 0.6, \gamma = 1, \delta = 0.5, \eta = 0.9,$ and $\kappa = 0.8$. This system was constructed for this study as a controlled test case exhibiting delayed cross-coupling and nontrivial oscillatory behavior; to our knowledge it does not correspond to a standard benchmark in the literature. The domain $[0,20]$ is divided into 20 uniform segments, {each trained with $N_f=1000$ interior points and $N_h=1500$ history samples.} Two fully connected sub-networks are used, one for $u(t)$ and one for $v(t)$. They are trained together using a single combined loss function, meaning that
the physics residuals for both equations contribute to the same objective
    $\mathcal{L}_{total} = \mathcal{L}_{hist} + 1.5\mathcal{L}_{DDE} + 2\mathcal{L}_{bc}.$
This ensures that $u(t)$ and $v(t)$ evolve consistently under the coupled DDE system. To illustrate the benefit of segmentation, we compare our model with three non-segmented Trial-PINN baselines of increasing capacity: PINN-1 has the same architecture as the segmented model, PINN-2 uses a $4 \times 20$ network with 5000 epochs, and PINN-3 uses a $5 \times 20$ network with 10000 epochs. 
\newline
Figure $\ref{fig:example3}$, shown over the interval [0,13], compares model predictions with the MATLAB {\texttt{dde23}} reference solution. The segmented method provides accurate amplitude and phase alignment for both $u(t)$ and $v(t)$ across the full domain, while the non-segmented models show noticeable phase lag and amplitude drift, particularly for $t<12.5$.
\begin{figure}[H]
  \centering
  \includegraphics[width=\columnwidth]{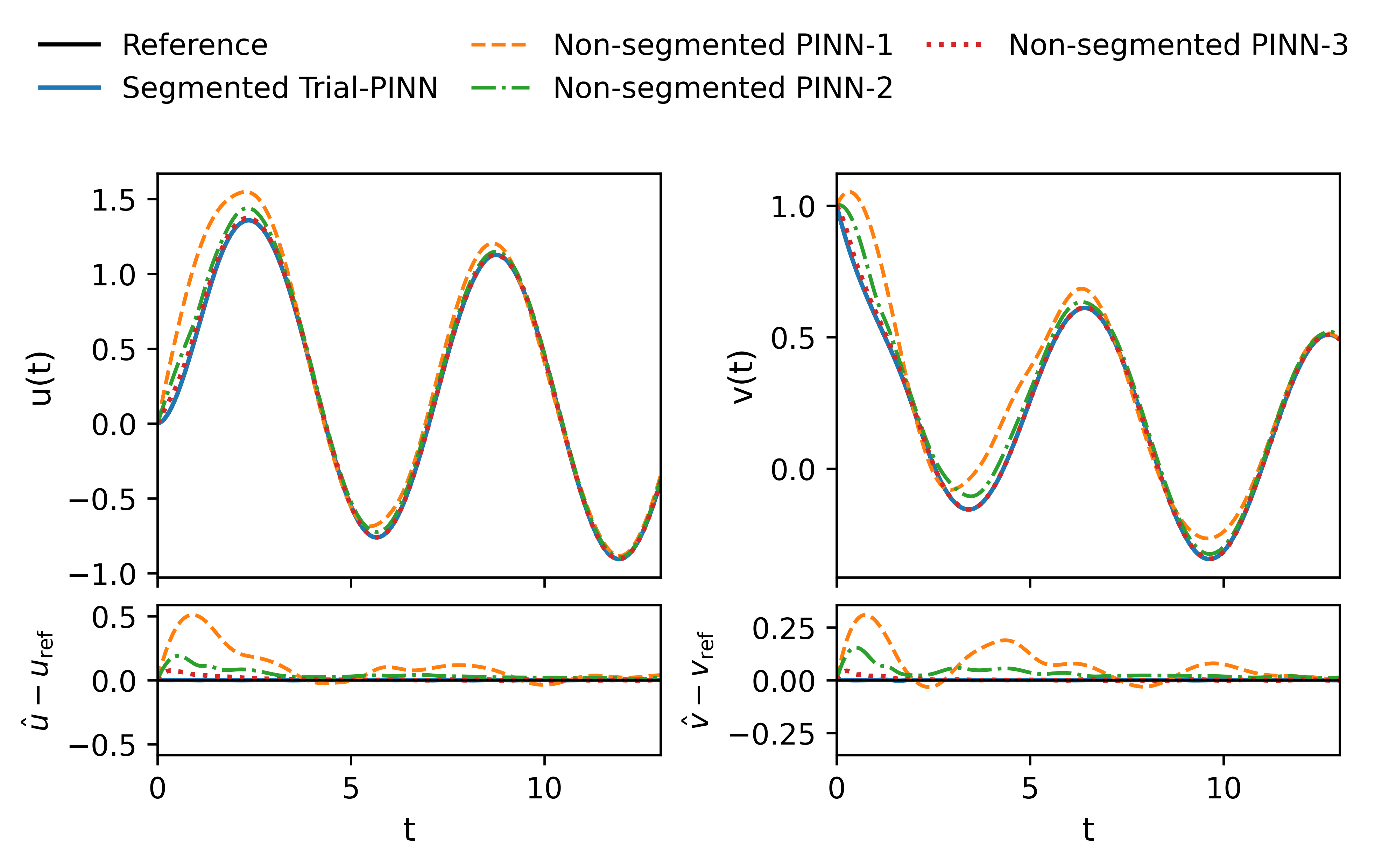}
  \caption{Comparison of segmented and non-segmented models for the coupled DDE \eqref{DDE-3}. Variable $u(t)$ in the left column, and $v(t)$ in the right column. 
  }
  \label{fig:example3}
\end{figure}
Table $\ref{tab:training_results_5}$ summarizes training cost and accuracy across all models. The segmented Trial-PINN achieves the lowest error overall: its RMSE $(6.6\times10^{-4})$ is nearly two orders of magnitude smaller than that of the deepest non-segmented model $(1.22\times10^{-2})$, while its MAE shows a similar improvement $(4.4\times10^{-7}$ vs. $1.49\times10^{-5})$. Although the deepest non-segmented PINN reduces error relative to smaller baselines, it still remains $\approx20\times$ less accurate in RMSE and $\approx 34\times$ less accurate in MAE, although it requires about half the training time of the segmented model (895 s vs. 450 s).
\begin{table}[h]
  \centering 
  \caption{Quantitative comparison of different models on the coupled DDE \eqref{DDE-3}}
  \label{tab:training_results_5}
  \begin{tabular}{p{3.9cm} p{1.15cm} p{.9cm}  p{0.8cm}}
    \hline
    Model & RMSE ($×10^{-3}$) & MAE ($×10^{-4}$) & Training Time (s)\\
    \hline
    Segmented Trial-PINN & 0.661 & 0.00044 & 449.92 \\
    Non-segmented Trial-PINN-1 & 140.45 & 197.25 & 118.24 \\
    Non-segmented Trial-PINN-2 & 81.49 & 66.41 & 144.60 \\
    Non-segmented Trial-PINN-3 & 12.21 & 1.49 & 895.46 \\
    \hline
  \end{tabular}
\end{table}
\end{example}

\section{Conclusion and Discussion}

In this paper, we introduced the segmented Trial-PINN framework for solving general nonlinear DDEs. The method integrates three key components: a differentiable history switch for handling delayed dependencies, a trial-solution formulation that enforces history constraints by construction, and a segmented collocation strategy that improves numerical stability and scalability. Together, these elements yield a 
physics-consistent approximation that remains continuous across subdomains while mitigating gradient instability from delayed feedback.

Three numerical examples illustrate where segmentation helps and where it does not. The first two examples are standard benchmarks in the DDE literature, while the coupled system provides a controlled test case; they provide standard testbeds for evaluating stability and delay-handling performance. In the harmonic DDE (Example ~\ref{ex:harmonic}), segmentation offers no accuracy benefit, showing that small, well-behaved systems are adequately handled by a standard non-segmented PINN. In contrast, for the nonlinear logistic DDE (Example ~\ref{ex:logistic}) and especially for the coupled system (Example ~\ref{ex:coupled}), the segmented Trial-PINN achieves substantially lower errors and remains more stable without requiring deeper or wider networks. Segmentation reduces local residuals, improves convergence in delayed regions, and provides robustness on large domains or multivariable systems. Although it introduces additional training cost due to solving multiple submodels, the overhead is modest relative to the improvements obtained on challenging delayed dynamics.

Methodologically, this work strengthens the generality of physics-informed neural networks by incorporating explicit delay modelling and structural continuity constraints. The resulting framework bridges traditional numerical solvers with data-driven learning and applies to a wide range of delayed dynamical systems in science and engineering. Future work will extend the framework to higher-order, time- and state-dependent delays, adaptive segmentation strategies, and hybrid PINN–numerical solver methods to further enhance efficiency and robustness. 

\end{document}